\RequirePackage{ifpdf}
\ifpdf 
\documentclass[pdftex]{sigma}
\else
\documentclass{sigma}
\fi

\numberwithin{equation}{section}

\begin{document}

\allowdisplaybreaks

\renewcommand{\thefootnote}{$\star$}

\renewcommand{\PaperNumber}{032}

\FirstPageHeading

\ShortArticleName{Vector Fields via Faber Polynomials}

\ArticleName{Vector Fields on the
Space of Functions Univalent\\ Inside the Unit Disk via Faber Polynomials\footnote{This paper is a
contribution to the Special Issue on Kac--Moody Algebras and Applications. The
full collection is available at
\href{http://www.emis.de/journals/SIGMA/Kac-Moody_algebras.html}{http://www.emis.de/journals/SIGMA/Kac-Moody{\_}algebras.html}}}

\Author{Helene AIRAULT}

\AuthorNameForHeading{H. Airault}

\Address{LAMFA CNRS UMR 6140, Insset, Universit\'e de Picardie Jules Verne,\\ 48
rue Raspail, 02100 Saint-Quentin (Aisne), France}

\Email{\href{mailto:helene.airault@insset.u-picardie.fr}{helene.airault@insset.u-picardie.fr}}

\ArticleDates{Received July 17, 2008, in f\/inal form March 07,
2009; Published online March 15, 2009}

\Abstract{We obtain the Kirillov vector f\/ields on the set of  functions $f$ univalent  inside the unit  disk, in terms of the Faber polynomials of $1/f(1/z)$. Our construction relies on the  generating  function for Faber polynomials.}


\Classification{17B66; 33C80; 35A30}

\section{Introduction}
The Virasoro algebra  has a representation
in the tangent bundle over the space of  functions univalent in the unit disk which are smooth on its boundary. This realization was obtained by A.A.~Kirillov
and D.V.~Yur'ev \cite{Kirillov,KirillovYuriev}  as f\/irst-order dif\/ferential
operators. Following \cite{Kirillov}, consider~$f(z)$  a
holomorphic  function univalent
in the unit disc $D=\{z\in \mathbb{C}  ; |z|\leq 1\}$, smooth up to the boundary of the disc and normalized by the conditions $f(0)=0$ and $f'(0)=1$, thus
\begin{gather}
f(z)=z\bigg(1+\sum_{n\geq 1}c_n z^n\bigg)\label{(1.1)}
\end{gather}
and the series \eqref{(1.1)} converges to $f(z)$ on $D$.  By De Branges theorem proving the Bieberbach conjecture, the coef\/f\/icients $(c_n)_{n\geq 1}$ lie in the inf\/inite-dimensional domain $|c_n|<n+1$, for $n\geq 1$.  In the following, we shall call this inf\/inite domain the manifold of coef\/f\/icients. On the other hand, we denote
\begin{gather*}
g(z)=b_0 z+ b_1+\frac{b_2}{ z}+\cdots +  \frac{b_n}{ z}+\cdots
\end{gather*}
 a  function univalent outside the unit disc. In order to study the representations of the Virasoro algebra \cite{Neretin}, A.A.~Kirillov  considered the action of vector f\/ields on the set of the dif\/feomorphisms of
the circle by perturbing the equation $f\circ \gamma=g$,
where $\gamma$ is a dif\/feomorphism of the circle. He obtained a sequence of vector f\/ields $L_{p}$, ($p$ being
positive or negative integer)  acting on the set of functions univalent  inside the unit disk. These vector f\/ields are expressed as
\begin{gather}
L_{-p}f(z)=\frac{f(z)^2}{ 2i\pi}\int_{\partial D}\frac{t^2f'(t)^2}{ f(t)^2}  \frac{1}{ f(t)-f(z)}  \frac{dt}{ t^{p+1}}=\phi_p(z)+z^{1-p}f'(z)\qquad\forall\, p\in \mathbb{Z},\label{(1.3)}
\end{gather}
where $z$ is inside the unit disk and the integral is a contour integral over the unit circle.
When $n>0$ is a positive integer, the action of $L_n$ is given by $L_nf(z)=z^{n+1}f'(z)$ since by evaluation of the contour integral, it is not dif\/f\/icult  to see that $\phi_{-n}(z)$
vanishes in~\eqref{(1.3)}. The term $\phi_p(z)$ comes from the residue at the pole $t=0$   and  $z^{1-p}f'(z)$ comes from the residue at $t=z$.
We have  $\phi_0(z)=-f(z)$
and $L_0f(z)=zf'(z)  -  f(z)$.  When $p>0$ the evaluation of the integral~\eqref{(1.3)} is more delicate since we have a non-vanishing residue at zero.
For $p\geq 0$, we f\/ind
\begin{gather*}
L_{-1}f(z)  =  f'(z)  -  1  -2c_1f(z),\nonumber\\
 L_{-2}f(z)  =  \frac{f'(z)}{z}  -  \frac{ 1}{f(z)}  -  3c_1  -  \big(4c_2-c_1^2\big)f(z),\nonumber\\
L_{-3}f(z)  =  \frac{f'(z)}{z^2}  -  \frac{1}{f(z)^2}  -  \frac{4c_1}{f(z)}  - \big(c_1^2+5c_2\big)  - (6c_3-2c_1c_2) f(z),
\end{gather*}
where the coef\/f\/icients $(c_j)_{j\geq 1}$ are the coef\/f\/icients of $f(z)$ given by~\eqref{(1.1)}. With the residue calculus, it has been made explicit in~\cite{Airault-Malliavin} that  $\phi_p(z)  = \Lambda_p(f(z))$ where $u\to \Lambda_p(u)$ is a function of the form
\begin{gather*}
\Lambda_p(u)=  -\frac{1}{u^{p-1}}  -\frac{(p+1)c_1}{u^{p-2}}  -  \cdots- a_p^p u. 
\end{gather*}
 The coef\/f\/icients of $\Lambda_p(u)$ depend on the $(c_j)_{j\geq 1}$ and have been calculated in \cite[Proposition~3.2]{Airault-Neretin}.  We note that $\phi_p(z)$ is obtained by eliminating the powers of $z^n$, $n\leq 1$ in $z^{1-p}f'(z)$; the elimi\-nation is done by expanding $z^{1-p}f'(z)$ in powers of $f(z)$, then substracting the part $\phi_p(z)$ of the series with powers  $f(z)^n$, $n\leq 1$. This
method is analogous to the  elimination of terms in power series  developed by Schif\/fer  for Faber polynomials~\cite{Schiffer}.
Let $h(z)$ be a univalent function holomorphic outside the unit disc except for a pole with residue equal to $1$ at
inf\/inity, thus for $|z|>1$,
\[
h(z)=z+b_1+b_2\frac{1}{z}+\cdots+b_n\frac{1}{z^{n-1}}+\cdots.
\]
Let $t\in \mathbb{C}$, at a neighborhood of $z=\infty$, we have the expansion
\[
\frac{zh'(z)}{h(z)-t}=\sum_{n=0}^{\infty}F_n(t)  z^{-n}.
\]
The function $F_n(t)$ is a polynomial of degree $n$ in the variable $t$ and is called the $n^{\rm th}$ Faber polynomial of the function $h$. Schif\/fer showed in~\cite{Schiffer} that the polynomial $F_n(t)$  is the unique polynomial in $t$ of degree $m$ such that
\[
F_m(h(z))=z^m+\sum_{n=1}^\infty a_{mn}z^{-n}.
\]
 The objective of this paper is to show that with a method analogous to that of M.~Schif\/fer for obtaining  Faber polynomials, we can recover Kirillov vector f\/ields $L_{-p}$ when $p\geq 0$. For this,
 let~$f(z)$ be a univalent function as in~\eqref{(1.1)}, for $p\geq 0$,  we start from $z^{1-p}f'(z)$, it expands in~$D$~as
\[
z^{1-p}f'(z)=z^{1-p}+  \sum_{1\leq n\leq p+ 1}n  c_n z^{n-p}  +  \hbox{terms in } z^n\qquad (n\geq 2).
\]
 In Section~\ref{section2},  the function $\Lambda_p(u)$ for $p\geq 0$, is constructed in such
a way that
\begin{gather}
z^{1-p}f'(z)+\Lambda_p(f(z))\label{(1.6)}
\end{gather}
expands in powers $z^n$ with $n\geq 2$. In fact, the function  $\Lambda_p(t)$ with respect to $z^{1-p}f'(z)$ plays the same role as the Faber polynomials  $F_m(t)$ with respect to $h(z)^m$. Then we prove that $z^{1-p}f'(z)  +  \Lambda_p(f(z))$ is equal to the expression~\eqref{(1.3)} of the vector f\/ield $L_{-p}f(z)$ found by Kirillov and Yur'ev.

Note that the method  of elimination of terms in power series as developed in  \cite{Schiffer} for Faber polynomials is a~formal calculation on series and does not require smoothness assumptions for~$f(z)$ at the boundary of the unit disk. Thus, it is conceivable to extend the calculations of the vector f\/ields for functions $f(z)$ which present a singularity at the boundary of $D$.
 This is our main motivation for adapting Schif\/fer's  method to the formal series $z^{1-p}f'(z)$.
  The regularity assumptions on $f(z)$ are stronger  in the case of the variational approaches developed in~\cite{Kirillov} or~\cite{Schaeffer-Spencer} since in the  variational case, it is assumed that $f(z)$ is smooth up to the boundary of the unit disc.

Let $\Lambda_p(u)$ as in~\eqref{(1.6)}, we  give an expression of $\Lambda_p(u)$, $p>0$,  in terms of the Faber polyno\-mials~$F_n(w)$ of the function $h(z)=  1/ f(1/z)$. We have
\[
\frac{ \xi h'(\xi)}{ h(\xi)-w}=
\sum_{n=0}^{\infty}F_n(w) \xi^{-n},
\] where $F_n(w)$ are the Faber polynomials associated to the function
$h$. In terms of $f$,
\begin{gather}
\frac{z f'(z)}{f(z)-w f(z)^2}=1+\sum_{n\geq 1}F_n(w)z^n,\label{(1.7)}\\
F_1(w)=w+c_1,  \qquad F_2(w)=w^2+ 2c_1 w+2c_2-c_1^2,\nonumber\\
F_3(w)=w^3+3c_1w^2 +3 c_2 w+c_1^3-3 c_1 c_2+3 c_3.\nonumber
\end{gather}
We f\/ind that the functions $u\to\Lambda_p(u)$
are determined by the
expansion
 \begin{gather}
 \frac{\xi^2 f'(\xi)^2}{ f(\xi)^2} \frac{u^2}{
 f(\xi)-u}  =  \sum_{p\geq 0}\Lambda_p(u)\xi^p\label{(1.8)}
 \end{gather}
  at a neighborhood of $\xi=0$, and  $u$ is a complex number.  For $p\geq 2$, we show that we can calculate the coef\/f\/icient $\Lambda_p(u)$  of  $\xi^p$ in the expansion~\eqref{(1.8)} as follows,
\begin{gather*}
\Lambda_p(u)  +  a_p^p  u  =  -  T_{p-1}\left(\frac{1}{u}\right),
\end{gather*}
where
\begin{gather*}
 T_{p-1}(w)= F_{p-1}(w)+2 c_1F_{p-2}(w)+3 c_2F_{p-3}(w)+\cdots +(p-1)c_{p-2}
F_{1}(w)+p c_{p-1}
\end{gather*} is determined by
\begin{gather*}
\frac{z f'(z)^2}{ f(z)-w f(z)^2}=1+\sum_{n\geq 1}T_n(w)z^n
\end{gather*}
and $a_p^p$,
\begin{gather}
\frac{z^2f'(z)^2}{f(z)^2}  =  1  +  \sum_{p\geq 1}a_p^p  z^p.\label{(1.12)}
\end{gather}
Then, we recover~\eqref{(1.3)} as Corollary, see~\eqref{(2.24)}.
In Section~\ref{section3}, we put
\begin{gather}
z^{1-p}f'(z)+\Lambda_p(f(z))=\sum_{n\geq 1}A_n^pz^{n+1}.\label{(1.13)}
\end{gather}
We prove that for any $u$ and $v$ in the unit disc, there holds
\begin{gather}
\sum_{k\geq 1}\sum_{p\geq 0}A_k^pu^pv^{k}  =  \frac{u^2f'(u)^2}{f(u)^2}  \frac{f(v)^2}{v[f(u)-f(v)]}  +  \frac{vf'(v)}{v-u}. \label{(1.14)}
\end{gather}
We say that the right hand side of~\eqref{(1.14)} is a generating function for the homogeneous polyno\-mials~$A_k^p$.
Note that the right hand side in~\eqref{(1.14)} has a meaning when~$u\to v$. All series obtained as expansions of a function are convergent inside their  disc of convergence which is determined by the singularities of the function. In Section~\ref{section4},  we identify as in~\cite{Kirillov} the vector f\/ields $(L_k)_{k\geq 1}$ with f\/irst order dif\/ferential operators on the manifold of coef\/f\/icients of functions univalent on $D$; as quoted before, this manifold comes from De Branges theorem, the coef\/f\/i\-cients~$(c_n)_{n\geq 1}$ lie in the inf\/inite-dimensional domain $|c_n|<n+1$, for $n\geq 1$. Some properties of this inf\/inite-dimensional manifold have been investigated in~\cite{Kirillov,KirillovYuriev}, for example K\"{a}hler structure or Ricci curvature. Here, we shall not develop the properties of this manifold. We only  examine the action of the  $(L_k)_{k\geq 1}$ on the functions  $\Lambda_p(u)$. The functions $\Lambda_p(u)$ have their coef\/f\/icients in this manifold. We f\/ind  that $L_k(\Lambda_{p+k}(u))  =  (2k+p)\Lambda_p(u)$ for $p\geq 1$.
 In  Section~\ref{section5}, we consider the reverse series $f^{-1}(z)$
  of $f(z)$, i.e.\ $f^{-1}\circ f(z)  =  z$. We prove that for $k>0$, $L_k[f^{-1}(z)]=  -  [f^{-1}(z)]^{k+1}$ and
\begin{gather*}
L_k\big[\big(f^{-1}(z)\big)^{-k}\big]  = k\qquad \text{if}\quad    k\geq 1,
\end{gather*}
  thus  the coef\/f\/icients in the expansion of $1/[f^{-1}(z)]^k$ in powers $z^n$, are  vectors $v(c_1, c_2,\dots)$ solution of  $L_k(v)=0$. On the other hand for $p\geq 1$, there holds
\begin{gather*}
L_{-p}\big[\big(f^{-1}(z)\big)^p\big]  = -p  -\Lambda_p(z)  \frac{d}{dz}\big [f^{-1}(z)\big]^p.
\end{gather*}

  \section[Elimination of terms in power series and the vector fields $L_{-p}f(z)$ in terms of the Faber polynomials of $1/f(1/z)$]{Elimination of terms in power series and the vector\\ f\/ields $\boldsymbol{L_{-p}f(z)}$ in terms of the Faber polynomials of $\boldsymbol{1/f(1/z)}$}\label{section2}

 Let $f(z)$ as in \eqref{(1.1)}, there exists a unique sequence of rational functions  $(\Lambda_p)_{p\geq 0}$ of the form
 \begin{gather}
 \Lambda_p(u)=\alpha_0 u  +  \alpha_1  + \frac {\alpha_2}{ u}  +  \cdots  -  \frac{1}{ u^{p-1}}, \label{(2.1)}
\end{gather}
  such that
\begin{gather}
z^{1-p}f'(z)  =  -  \Lambda_p[f(z)]  + \mbox{series of  terms  in} \  z^k, \qquad  k\geq2. \label{(2.2)}
\end{gather}
 To prove the existence of $\Lambda_p(u)$,  we expand $z^{1-p}f'(z)$ in powers of $f(z)$. Then $-  \Lambda_p(f(z))$ is the sum of terms with powers of $f(z)^n$ such that $n\leq 1$. The unicity of the function $\Lambda_p(u)$ satisfying \eqref{(2.1)}, \eqref{(2.2)} results from matching equal powers of $z$ in the expansions  \eqref{(2.1)}, \eqref{(2.2)}.  We can calculate directly
\begin{gather*}
\Lambda_0(w)=-w, \qquad \Lambda_1(w)=-1-2 c_1 w , \qquad \Lambda_2(w)=-\frac{1}{w}-3 c_1-
\big(4c_2-c_1^2\big) w,\\
\Lambda_3(w)=-\frac{1}{w^2}-4 c_1 \frac{1}{w}-\big(c_1^2+5 c_2\big) - (6c_3-2c_1c_2) w, \qquad \ldots.
\end{gather*}
For $p=0$, we have
\[
f(z)=z+c_1z^2+\cdots, \qquad  zf'(z)=z+2c_1z^2+3c_2z^3+\cdots,
\]
 thus {\samepage
\[
zf'(z)-f(z)=c_1z^2+2c_2z^3+\cdots=L_0[f(z)]
\]
 and $\Lambda_0(u)=  -  u$.}

For $p=1$,
\[
f'(z)=1+2c_1z+3c_2z^2+\cdots,
\]
 then
 \[
 f'(z)-1-2c_1f(z)=\big(3c_2-2c_1^2\big)z^2+\cdots  =  L_{-1}[f(z)].
 \]
 We obtain $\Lambda_1(u)=-1-2c_1u$.

For $p=2$,
\[
\frac{ f'(z)}{ z}=  \frac{1}{z}  +2c_1  +3c_2z +4c_3z^2+\cdots,
\] then
\begin{gather*}
\frac{ f'(z)}{ z}-\frac{ 1}{ f(z)} =3c_1  +(3c_2-G_2)z+\cdots\qquad \mbox{with} \quad   G_2=c_1^2-c_2,\\
 \frac{f'(z)}{z}-  \frac{ 1}{ f(z)} - 3c_1  -(3c_2-G_2)f(z)  =  \mbox{coef\/f\/icients} \cdot z^2+  \cdots.
\end{gather*}
This gives $ \Lambda_2(u)=  -  \frac{ 1}{ u}  -  3c_1  -  (3c_2-G_2)u$.
In the following theorem, we prove~\eqref{(1.8)}.

\begin{theorem}\label{theorem2.1}  $\Lambda_p(w)$
is expressed in terms of the Faber polynomials of $h(z)=\frac{ 1}{f(1/z)}$; we have
\begin{gather}
\Lambda_p(u)  +  a_p^p  u  =  -  T_{p-1}\left(\frac{1}{u}\right),\label{(i)}\\
T_{p-1}(w)= F_{p-1}(w)+2 c_1F_{p-2}(w)+3 c_2F_{p-3}(w)+\cdots +(p-1)c_{p-2}
F_{1}(w)+p c_{p-1},\label{(ii)}\\
 \frac{z f'(z)^2}{ f(z)-w f(z)^2}=1+\sum_{n\geq 1}T_n(w)z^n, \label{(iii)}\\
 \frac{z^2 f'(z)^2}{ f(z)^2}=1+\sum_{p\geq 1} a_p^p z^p. \label{(iv)}
\end{gather}
The functions  $(\Lambda_p(w))_{p\geq 0}$ are
given by
\begin{gather}
\frac{\xi^2 f'(\xi)^2}{ f(\xi)^2} \frac{u^2}{
 f(\xi)-u}  =  \sum_{p\geq 0}\Lambda_p(u)\xi^p. \label{(v)}
\end{gather}
\end{theorem}

\begin{proof}
Consider the function
\begin{gather*}
h(z)=\frac{1}{ f(\frac{ 1}{
 z})}=  z-c_1+\big(c_1^2- c_2\big)\frac{1}{ z}
+\big(2 c_1c_2-c_3-c_1^3\big)\frac{1}{ z^2} \nonumber\\
\phantom{h(z)=}{} +\big(2 c_1 c_3-c_4+c_2^2-3 c_1^2c_2+ c_1^4\big)\frac{1}{z^3}+\cdots.
\end{gather*}
As in \cite{Schiffer},  we have
\[
\frac{\xi h'(\xi)}{ h(\xi)-w}=
\sum_{n=0}^{\infty}F_n(w) \xi^{-n},
\] where $F_n(w)$ are the Faber polynomials associated to the function~$h$. Since $h(z)=1/f(1/z)$, we have~\eqref{(1.7)}.
If we take the derivative of~\eqref{(1.7)} with respect to $w$ and then
integrate with respect to $z$, we obtain
\begin{gather*}
\frac{f(z)}{(1-wf(z))}=\sum_{n\geq 1}F'_n(w)\frac{z^n}{ n}. 
\end{gather*}
Moreover
\begin{gather*}
F_n(h(z))=z^n+\sum_{k=1}^{\infty}\beta_{n,k}z^{-k},
\end{gather*}
where the $\beta_{n,k}$ are the Grunsky coef\/f\/icients of $h$, see~\cite{Schiffer}. In terms of $f(z)$,
\begin{gather}
K(u,v)=\log \frac{\frac{ 1}{f(u)}  -  \frac{ 1}{ f(v)}}{ \frac{ 1}{ u}-\frac{ 1}{ v}}
=-\sum_{n\geq 1}\sum_{k\geq 1}\frac{1}{ n}\beta_{n,k}u^nv^k.
\label{(2.6)}
\end{gather}
Because of the symmetry in $(u,v)$ of the left hand side in \eqref{(2.6)}, we see that
\begin{gather}
\frac{1}{ n}\beta_{n,k}  =  \frac{1}{ k}\beta_{k,n}. \label{(2.7)}
\end{gather}
 Thus for $n\geq 1$,
\begin{gather}
F_n\left( \frac{ 1}{ f(z)}\right) =z^{-n}+\sum_{k=1}^{\infty}\beta_{n,k}z^{k}.\label{(2.8)}
\end{gather}
 We rewrite \eqref{(2.8)} as
 \[
 z^{-n}=F_n\left(\frac{ 1}{f(z)}\right)-
\sum_{k=1}^{\infty}\beta_{n,k}z^{k}.
\]
On the other hand, if $p> 1$,
\begin{gather}
z^{1-p}f'(z)=z^{1-p}(1+2 c_1 z+ 3c_2 z^2+\cdots +(n+1)c_nz^n
+\cdots\nonumber\\
\phantom{z^{1-p}f'(z)}{} =\frac{1}{ z^{p-1}}+\frac{2c_1}{ z^{p-2}}+\frac{3 c_2}{ z^{p-3}}+\cdots
+\frac{(p-1)c_{p-2}}{ z } \nonumber\\
\phantom{z^{1-p}f'(z)=}{}+ pc_{p-1}+(p+1)c_p z
+  \sum_{k\geq 1}(p+k+1)c_{p+k}z^{k+1}.\label{(2.9)}
\end{gather}
We replace in~\eqref{(2.9)} the negative powers of $z$ by their expressions given in~\eqref{(2.8)}. We obtain
\begin{gather}
z^{1-p}f'(z)=  F_{p-1}\left(\frac{1}{f(z)}\right)
+2 c_1F_{p-2}\left(\frac{1}{f(z)}\right)
+  3 c_2F_{p-3}\left(\frac
{1}{ f(z)}\right)+\cdots \nonumber\\
\phantom{z^{1-p}f'(z)=}{} +(p-1)c_{p-2} F_{1}\left(\frac{1}{f(z)}\right)
 +p c_{p-1}+(p+1)c_p z\nonumber\\
 \phantom{z^{1-p}f'(z)=}{}
-\big[ \beta_{p-1,1}+ 2c_1 \beta_{p-2,1}+\cdots
+(p-1)c_{p-2}\beta_{1,1}\big]  z \nonumber\\
\phantom{z^{1-p}f'(z)=}{}
 +\sum_{k\geq 1}\big[ (p+ k+1)c_{p+k}-\big [ \beta_{p-1,k+1}+ 2c_1 \beta_{p-2,k+1}+\cdots\nonumber\\
\phantom{z^{1-p}f'(z)=}{}+(p-1)c_{p-2}\beta_{1,k+1}\big] \big]  z^{k+1}.\label{(2.10)}
\end{gather}
For $p\geq 2$, we consider
\begin{gather*}
T_{p-1}(w)= F_{p-1}(w)+2 c_1F_{p-2}(w)+3 c_2F_{p-3}(w)+\cdots +(p-1)c_{p-2}
F_{1}(w)+p c_{p-1}. 
\end{gather*}
From the expansion of $f'(z)$ and \eqref{(1.7)}, we obtain
\begin{gather*}
\frac{z f'(z)^2}{ f(z)-w f(z)^2}=1+\sum_{n\geq 1}T_n(w)z^n,\\ 
T_0(w)=1, \qquad T_1(w)=w+3c_1, \qquad T_2(w)=w^2+4 c_1w+\big(c_1^2+5 c_2\big),\nonumber\\
T_3(w)=w^3+5 c_1w^2+\big(4c_1^2+6 c_2\big)w  -   c_1^3+4c_1c_2+7c_3, \qquad \ldots.\nonumber
\end{gather*}
We write \eqref{(2.10)} as
\begin{gather}
z^{1-p}f'(z)-T_{p-1}\left( \frac{1}{ f(z)}\right) =
\sum_{k\geq 1}B_{k-1}^p z^k \nonumber\\
\qquad{} =(p+1)c_pz  -\big [ \beta_{p-1,1}+ 2c_1 \beta_{p-2,1}+\cdots
+(p-1)c_{p-2}\beta_{1,1}\big ]  z+\sum_{k\geq 1}B_k^p z^{k+1}
\label{(2.13)}
\end{gather}
with
\begin{gather}
B_k^p=(p+k+1)c_{p+k}-\big [ \beta_{p-1,k+1}+ 2c_1 \beta_{p-2,k+1}
+\cdots +(p-1)c_{p-2}\beta_{1,k+1}\big ] .\label{(2.14)}
\end{gather}
At this step we have eliminated all the negative powers and the constant term in such a way that the series on the right side of~\eqref{(2.13)} has only terms in $z^n$ with $n\geq 1$. To eliminate the term in $z$ in order to have only terms in $z^n$, $n\geq 2$, we put for $p\geq 1$,
\begin{gather}
z^{1-p}f'(z)-T_{p-1}\left( \frac{1}{ f(z)}\right) -a_p^p f(z)
=\sum_{k\geq 1}A_k^p z^{k+1}. \label{(2.15)}
\end{gather}
From \eqref{(2.13)}, we see that the coef\/f\/icients of $z$,  $a_p^p$ are determined with
 $a_1^1=2c_1$ and if $p>1$,
\begin{gather}
a_p^p=-\big [ \beta_{p-1,1}+ 2c_1 \beta_{p-2,1}+\cdots
+(p-1)c_{p-2}\beta_{1,1}\big ] +(p+1)c_p. \label{(2.16)}
\end{gather}
For $p\geq 1$, $k\geq 1$, we put
\begin{gather}
A_k^p=B_k^p-a_p^pc_k. \label{(2.17)}
\end{gather}
With the convention $c_0=1$, $B_0^p=a_p^p$, we have $A_0^p=0$.
Now, we prove that  $(a_p^p)$
is given by~\eqref{(1.12)} or \eqref{(iv)}. For this, we consider \eqref{(2.8)} with $n=1$, it gives
\begin{gather}
F_1\left(\frac{1}{ f(z)}\right)  =  \frac{1}{ z}  +  \sum_{k\geq 1}\beta_{1,k}z^k  =  \frac{1}{ f(z)}  +  c_1. \label{(2.18)}
\end{gather}
Because of the symmetry \eqref{(2.7)} in the Grunsky coef\/f\/icients, taking the derivative with respect to $z$ in~\eqref{(2.18)}, we obtain,
\begin{gather}
-  \frac{1}{ z^2}  +\sum_{k\geq 1}\beta_{k, 1}z^{k-1}  =  - \frac {f'(z)}{f(z)^2}. \label{(2.19)}
\end{gather}
Then we multiply the power series $f'(z)=1+2c_1z+3c_2z^2+\cdots$ by the power series in \eqref{(2.19)}, it gives  \eqref{(2.16)} on one side and \eqref{(1.12)} or \eqref{(iv)} on the other side.
To prove \eqref{(1.8)} or \eqref{(v)}, we put
\begin{gather*}
M_p(w)=-T_{p-1}(w)-a_p^p  \frac
{ 1}{ w}\qquad  \mbox{for}\quad   p\geq 1\qquad \mbox{and} \qquad
M_0(w)= - \frac{ 1}{ w}, 
\end{gather*} then with \eqref{(iii)} and \eqref{(iv)},
\begin{gather*}
\sum_{p\geq 0}M_p(w)\xi^p=-\frac{\xi^2\,f'(\xi)^2}{ f(\xi)^2}  \frac{1}{
w  (1-w f(\xi))}
\end{gather*} and from \eqref{(2.15)}, we see that
\begin{gather*}
z^{1-p}f'(z)+M_p\left(\frac{1}{f(z)}\right)=\sum_{n=1}^{\infty}A_n^p z^{n+1}.
\end{gather*}
We put $\Lambda_p(w)=M_p(\frac{1}{w})$. We obtain \eqref{(v)} or \eqref{(1.8)}.
\end{proof}

\begin{corollary}\label{corollary 2.2}
Let
$\phi_p(z)=\Lambda_p(f(z))$, then
\begin{gather}
\sum_{p\geq 0}\phi_p(z)\xi^p=\frac{
 \xi^2 f'(\xi)^2}{ f(\xi)^2}  \frac{ f(z)^2}{ (f(\xi)- f(z))}\label{(2.23)}
 \end{gather}
 and
\begin{gather}
\frac{f(z)^2}{ 2 i\pi}\int \frac{\xi^2\,f'(\xi)^2}{ f(\xi)^2}  \frac{1}{(f(\xi)- f(z))}  \frac{d\xi}{\xi^{p+1}}=\phi_p(z)+ z^{1-p}f'(z).
\label{(2.24)}
\end{gather}
\end{corollary}

\begin{proof} \eqref{(2.23)} is the immediate consequence of \eqref{(1.8)}.
 To prove \eqref{(2.24)}, we see that $\phi_p(z)$ is the coef\/f\/icient of $\xi^p$
in the expansion of
$\frac{\xi^2 f'(\xi)^2}{f(\xi)^2} \frac{f(z)^2}{(f(\xi)- f(z))}$ in powers of $\xi$.
We calculate the contour integral in \eqref{(2.24)} with the residue method. The term $z^{1-p}f'(z)$ comes from the residue at $\xi=z$.
\end{proof}

\section[The generating function for the $A_n^p$]{The generating function for the $\boldsymbol{A_n^p}$}\label{section3}

We put, see \eqref{(1.13)},
\[
z^{1-p}f'(z)+\Lambda_p(f(z))=\sum_{n\geq 1}A_n^pz^{n+1}.
 \]
 With \eqref{(2.17)}, \eqref{(2.14)} and \eqref{(1.12)}, we have obtained $A_n^p$ explicitly in terms of the Grunsky coef\/f\/i\-cients~$\beta_{n,k}$ of $h(z)=1/f(1/z)$ and in terms of the coef\/f\/icients $(c_j)_{j\geq 1}$ of $f(z)$.
In this section, we prove that $A_n^p$ are given by \eqref{(1.14)}. We consider for $|\xi|<|z|$ the series
\begin{gather}
\sum_{p\geq 0}z^{1-p}f'(z)\xi^p  =  zf'(z) \sum_{p\geq0}\left(\frac{\xi}{z}\right)^p  =  \frac{z^2 f'(z)}{z-\xi}\label{(3.1)}
\end{gather}
and, see~\eqref{(1.8)},
\begin{gather}
\sum_{p\geq 0}\Lambda_p(f(z))\xi^p  =  \frac{\xi^2f'(\xi)^2}{f(\xi)^2}  \frac{f(z)^2}{f(\xi)-f(z)}.\label{(3.2)}
\end{gather}
Adding \eqref{(3.1)} and \eqref{(3.2)}, then dividing by $z$, we f\/ind  for  $|\xi|<|z|$,
\begin{gather*}
\sum_{p\geq 0}\sum_{k\geq 1}A_k^pz^k\xi^p  =  \frac{\xi^2f'(\xi)^2}{f(\xi)^2}  \frac{f(z)^2}{z(f(\xi)-f(z))}  +  \frac{z f'(z)}{z-\xi},
\end{gather*}
which is \eqref{(1.14)} for $|u|<|v|$.
Below, we prove that \eqref{(1.14)} is true for any $u$ and $v$ in the unit disk.

\begin{theorem}\label{theorem 3.1}
 The polynomials $A_k^p$ defined by \eqref{(1.13)} satisfy  \eqref{(1.14)}.
\end{theorem}
\begin{proof}
Taking the derivative of \eqref{(2.6)} with respect to $u$ yields
\begin{gather}
\frac{u f'(u)}{ f(u)}  \frac{f(v)}{ f(v)-f(u)} -  \frac{v}{ v-u}=\sum_{n\geq 1}
\sum_{k\geq 1}\beta_{n,k}u^nv^k.\label{(3.4)}
\end{gather}
Multiplying \eqref{(3.4)} by $f'(u)$, we deduce that
\begin{gather}
\frac{u f'(u)^2}{ f(u)}  \frac{f(v)}{f(v)-f(u)}  -  \frac{vf'(u)}{ v-u}=
\sum_{n\geq 1}\sum_{k\geq 1}\big [ \beta_{n,k}+2 c_1 \beta_{n-1,k}+
\cdots +n c_{n-1}\beta_{1,k}\big ] u^nv^k. \label{(3.5)}
\end{gather}
Consider the homogeneous polynomials
 $B_k^p$  given by \eqref{(2.14)} for $p>1$ and $B_k^1=(k+2)c_{k+1}$. We
rewrite \eqref{(3.5)} as
\begin{gather}
\sum_{k\geq 0}\sum_{p\geq 2}\big(B_k^p-(p+k+1) c_{p+k}\big)u^{p-1}v^{k+1}=  -  \frac
{u f'(u)^2}{ f(u)}\frac{f(v)}{ f(v)-f(u)}
+  \frac{vf'(u)}{ v-u}. \label{(3.6)}
\end{gather}
Moreover, since $B_k^1=(k+2)c_{k+1}$, we can write the sum in \eqref{(3.6)}, starting from
$p=1$. On the other hand,
\[
-
\frac{ vf'(u)}{ v-u}  +  \frac{  vf'(v)}{  v-u} =  \frac{v}{ v-u}(f'(v)-f'(u))=v\sum_{k\geq 1}(k+1)c_k  \frac{v^k-u^k}{ v-u}
\]
thus
\begin{gather*}
\sum_{k\geq 0}\sum_{p\geq 1}(p+k+1) c_{p+k}u^{p-1}v^{k+1}=  -  \frac
{ vf'(u)}{ v-u}  +  \frac{ vf'(v)}{ v-u}.
\end{gather*}
 Finally, we obtain
\begin{gather*}
\sum_{k\geq 0}\sum_{p\geq 1}B_k^p u^{p}v^{k+1}=  -
\frac{u^2 f'(u)^2}{f(u)}  \frac{f(v)}{ f(v)-f(u)}  +  \frac{uvf'(v)}{ v-u}.
\end{gather*}
Since $A_k^p=B_k^p-a_p^pc_k$ as in \eqref{(2.17)}, $A_0^p=0$ and since \eqref{(iv)} is true,  we have
\begin{gather*}
\sum_{k\geq 1}\sum_{p\geq 1}A_k^p u^{p}v^{k+1}=  \frac{u^2 f'(u)^2}{
f(u)^2}  \frac{f(v)^2}{ f(u)-f(v)}+f(v)+  \frac{uv}{ v-u}f'(v).
\end{gather*}
We divide by $v$. Since
\[
\sum_{k\geq 1}\sum_{p\geq 0}A_k^p u^{p}v^{k}=
\sum_{k\geq 1}\sum_{p\geq 1}A_k^p u^{p}v^{k}+\sum_{k\geq 1}
A_k^0 v^{k}
\]
and
\[
\sum_{k\geq 1}A_k^0 v^{k}=\sum_{k\geq 0}k c_k v^k=f'(v)-  \frac{ f(v)}{v},
\] we obtain \eqref{(1.14)}.
\end{proof}

\section[The differential operators  $(L_k)_{k\geq 1}$ on the  functions $\Lambda_p(u)$]{The dif\/ferential operators  $\boldsymbol{(L_k)_{k\geq 1}}$ on the  functions $\boldsymbol{\Lambda_p(u)}$}\label{section4}

We identify the set of  functions univalent on the unit disk with the set of their coef\/f\/icients
 via the map
\begin{gather*}
f(z)=z\bigg( 1+\sum_{n\geq 1}c_n z^n\bigg) \quad \rightarrow \quad (c_1, c_2,
\dots, c_n,\dots).
\end{gather*}
For $k\geq 1$, we put $\partial_k=  \frac{\partial}{\partial c_k}$. Following~\cite{Kirillov}, we consider the partial dif\/ferential operators
\begin{gather}
L_k=\partial_k+\sum_{n=1}^{\infty}(n+1)c_n\partial_{n+k}.
\label{(4.2)}
\end{gather}
We have $z^{1+k}f'(z)=L_k[f(z)]$ and $\partial_n L_j=L_j
\partial_n+(n+1)\partial_{n+j}$.
We put
\begin{gather*}
L_0=\sum_{n\geq 1}n c_n\partial_n\qquad \mbox{and  for}\quad
k\geq 1\qquad L_{-k}=\sum_{n\geq 1} A_n^k \partial_n . 
\end{gather*}
The operators $\partial_k$ are related to the $(L_{k+p})_{p\geq 0}$ as follows, see \cite{Airault-Ren},
\begin{lemma}\label{lemma 4.1} For $k\geq 1$,
\begin{gather}
\partial_k=L_k -2c_1L_{k+1}+\big(4c_1^2-3c_2\big)L_{k+2}+\cdots + B_n L_{k+n}
+\cdots,\label{(4.4)}
\end{gather}
where the $(B_n)_{n\geq 1}$ are independent of $k$.
We calculate the $B_n$, $n\geq0$, with
\begin{gather*}
\frac{1}{ f'(z)}=1+\sum_{n\geq 1}B_n z^n.
\end{gather*}
\end{lemma}
\begin{proof}
We  verify \eqref{(4.4)} on $f(z)$. We have $L_k[f(z)]=z^{k+1}f'(z)$.
Since
\begin{gather*}
\partial_k [f(z)]=z^{k+1}\qquad \mbox{and}\qquad \partial_k[f'(z)]
=(k+1)z^k, 
\end{gather*}
we have to prove
\[
z^{k+1}=z^{k+1}f'(z)-2c_1 z^{k+2}f'(z)+\cdots +
B_n z^{k+n} f'(z).
\] We divide by $z^{k+1}$ and we obtain \eqref{(4.4)}.
\end{proof}

By considering expansions as in \cite{Airault-Malliavin},  we obtain with \eqref{(1.8)}, the action of $(L_k)_{k\geq 1}$
on the  functions $(\Lambda_p(u))_{p\geq 0}$.
\begin{theorem}\label{theorem 4.2}
Let $\Lambda_p(u)$ as in \eqref{(1.8)}, then
$L_k(\Lambda_n(u))=0$ if $1\leq n<k$, $\, \, L_k(\Lambda_k(u))=\, -\, 2ku$ for $k\geq 1$,
\begin{gather}
L_k(\Lambda_{p+k}(u))=(2k+p)\Lambda_p(u)\qquad\mbox{for}\quad p\geq 0\quad  \mbox{and} \quad k\geq 1.\label{(4.7)}
\end{gather}
\end{theorem}

\noindent
{\bf Rremark.}
The functions $\Lambda_p(u)$ are of the form (see also \cite[(A.1.7)]{Airault-Malliavin}  and \cite[Proposition~3.2]{Airault-Neretin})
\begin{gather*}
\Lambda_p(u)=  -  \bigg[  \frac{1}{ u^{p-1}}   + (p+1)c_1\frac{1}{ u^{p-2}}  +  \cdots  + \frac{ (p+n)c_n+\gamma_n(p)}{ u^{p-n-1}}  +  \cdots   + (2pc_p +\gamma_p(p))  u\bigg],
\end{gather*}
where for $2\leq n\leq p$, $\gamma_n(p)(c_1, c_2, \dots, c_{n-1})$ are homogeneous polynomials of degree $n$ in the variables $(c_1, c_2,\dots, c_{n-1})$ with $c_j$ having weight $j$.
With $k=1$ in \eqref{(4.7)}, we have
\begin{gather}
\left[\frac{\partial}{\partial c_1}+2c_1\frac{\partial}{\partial c_2}+  3c_2\frac{\partial}{\partial c_3}+\cdots\right](\Lambda_{p+1}(u))=(p+2)\Lambda_p(u). \label{(4.9)}
\end{gather}
Since $\Lambda_0(u)=-u$, one can calculate the sequence $(\Lambda_p(u))_{p\geq 1}$ recursively by identifying equal powers of $u$ in~\eqref{(4.9)}.

\section[The  $(L_k)_{k\geq 1}$ and the reverse series of $f(z)$]{The  $\boldsymbol{(L_k)_{k\geq 1}}$ and the reverse series of $\boldsymbol{f(z)}$}\label{section5}

As in Section~\ref{section4}, we consider the dif\/ferential operators
$(L_k)_{k\geq 1}$ given by~\eqref{(4.2)}.
We denote   $f^{-1}(z)$  the inverse function of $f(z)$, ($f^{-1}\circ f  ={\rm Identity})$, we say also reverse series of $f(z)$.
For any integer $q$, consider the series
\[
\big(f^{-1}(z)\big)^q=z^q\bigg(1 +\sum_{n\geq 1} \delta_n^q z^n\bigg),
\]
where $\delta_n^q$ are homogeneous polynomials in the variables $(c_1,c_2,\dots, c_n,\dots)$, the coef\/f\/icients of~$f(z)$. Then
\[
L_p\big[\big(f^{-1}(z)\big)^q\big]  =  z^q \sum_{n\geq 1} L_p\big [\delta_n^q\big ] z^n.
\]

\begin{theorem}\label{theorem 5.1}
Let $f^{-1}(z)$ be the reverse series of $f(z)$,  then
\begin{gather}
L_k\big[f^{-1}(z)\big]=  -  \big[f^{-1}(z)\big]^{k+1}\qquad  \mbox{for}\quad  k\geq 1,\label{(i)+}\\
L_0\big[f^{-1}(z)\big]= - f^{-1}(z) + z \big(f^{-1}\big)'(z),\nonumber\\
 L_{-1}\big[f^{-1}(z)\big]= -1  +  (1+2c_1z) \big( f^{-1}\big)'(z),\nonumber\\
L_{-p}\big[f^{-1}(z)\big]= - \big[f^{-1}(z)\big]^{1-p}  -  \Lambda_p(z) \big(f^{-1}\big)'(z)\qquad\mbox{for}\quad  p\geq 2. \label{(ii)+}
\end{gather}
In particular there exists a unique rational function of $z$, which is $\Lambda_p(z)$, such that $[f^{-1}(z)]^{1-p}  +  \Lambda_p(z) (f^{-1})'(z)$
expands in a Taylor series $\sum\limits_{n\geq 2}a_nz^n$ with powers $z^n$, $n\geq 2$.
Moreover
\begin{gather}
L_k\big[\big(f^{-1}(z)\big)^{-k}\big]  = k\qquad  \mbox{for}\quad  k\geq 1,\label{(iii)+}\\
L_{-p}\big[\big(f^{-1}(z)\big)^p\big] = -p  -\Lambda_p(z)\times \frac{d}{dz}\big[f^{-1}(z)\big]^p\qquad  \mbox{for}\quad  p\geq 1. \label{(iv)+}
\end{gather}
\end{theorem}

\begin{proof} $f\circ f^{-1}(z)=z$, dif\/ferentiating with
 the vector f\/ield $L_k$,
 \begin{gather*}
 (L_kf)\big(f^{-1}(z)\big)  +  f'\big(f^{-1}(z)\big) L_k\big[f^{-1}(z)\big] =  0. 
 \end{gather*}
 For $k\geq 1$,
 \[
 L_k[f(z)]=z^{1+k}f'(z),
 \] thus
 \[
 L_kf\big(f^{-1}(z)\big)=f'\big(f^{-1}(z)\big)\big[f^{-1}(z)\big]^{1+k}
 \]
  and
 this gives \eqref{(i)+}. Then we obtain \eqref{(iii)+} because
\[
L_k\big[\big(f^{-1}(z)\big)^{-k}\big] =- k f^{-1}(z)^{-k-1}  L_k\big[f^{-1}(z)\big].
 \]
Similarly,
\[
L_{-p}f(z)=z^{1-p}f'(z)+\Lambda_p(f(z))
\] for $p\geq 0$, we f\/ind \eqref{(ii)+} and we deduce \eqref{(iv)+} from
\begin{gather*}
L_{-k}\big[\big(f^{-1}(z)\big)^p\big] =  - p\big(f^{-1}(z)\big)^{p-k} - \Lambda_k(z)  \frac{d}{dz}\big(f^{-1}(z)\big)^p.\tag*{\qed} 
\end{gather*}\renewcommand{\qed}{}
 \end{proof}

\pdfbookmark[1]{References}{ref}
\LastPageEnding

\end{document}